# Game theory with integral equations as state dynamics

S. A. Belbas, Mathematics Dept., University of Alabama, Tuscaloosa, Al. 35487-0350, USA.
*e-mail*: SBELBAS@GMAIL.COM

*Abstract.* We formulate and analyze game-theoretic problems for systems governed by integral equations. For Volterra integral equations, we obtain and prove necessary and sufficient conditions for linear-quadratic problems, and for problems that are linear-quadratic in the control. Also, we obtain necessary conditions for one type of pursuit-evasion Volterra games.

1. <u>Introduction.</u>

The theory of differential games is a well-developed area, within the broader context of systems theory and control theory.

The present paper deals with game theory for systems described by integral equations. The results and techniques developed herein also apply to differential games, since Volterra integral equations contain also systems governed by ordinary differential equations. They also contain a variant of the particular type of pursuit and evasion policies formulated by Pontryagin in [P2]. Pontryagin defined a policy for one player as a "functional" of the recent history of the trajectory and the actions of the other player; it is reasonable to specialize this "functional" into a Volterra integral operator. When the evolution of the state is governed by a Volterra integral equation, an interpretation of each player's policy as a causal operator of the entire history of the game is already embedded into the model, and, in that sense, it extends Pontryagin's formulation.

We mention here that the study of differential games on the basis of necessary conditions of Pontryagin's extremum principle (usually called "maximum principle" in control theory) was initiated in Pontryagin's paper [P1]. Although today differential games are studied mostly by dynamic programming methods, the study of games for integral equations again requires necessary conditions and extremum principles, since integral equations are not well-suited for dynamic programming.

We note that a variety of problems in Economics, Population Dynamics, and related areas, have been modelled via integral equations. For such models, the related game-theoretic problems are naturally games with integral equations as state dynamics.

In this paper, we deal mostly with linear-quadratic (LQ) integral games and games that are linear-quadratic n the controls (LQC) but generally nonlinear in the state. We note that LQ differential games have been well-studied in the research literature, and merely as a sample we mention [B, E]. The theory of necessary conditions for optimal control problems for integral equations is well established, for example [S], and we shall make freely use of optimal control results in our study of game theory for integral equations.

## 2. Quadratic forms over $L^2$ spaces.

We consider the following quadratic functional over $L^2(G, R^m \times R^n)$, where $G$ is a bounded open set in $d$ – dimensional Euclidean space:

$$E(w_1, w_2) := \tfrac{1}{2} \int_G [w_1^T(x) \; w_2^T(x)] \begin{bmatrix} K_{11}(x) & K_{12}(x) \\ K_{21}(x) & K_{22}(x) \end{bmatrix} \begin{bmatrix} w_1(x) \\ w_2(x) \end{bmatrix} dx +$$

$$+ \tfrac{1}{2} \int_G \int_G [w_1^T(x) \; w_2^T(x)] \begin{bmatrix} L_{11}(x,y) & L_{12}(x,y) \\ L_{21}(x,y) & L_{22}(x,y) \end{bmatrix} \begin{bmatrix} w_1(y) \\ w_2(y) \end{bmatrix} dx\, dy + \int_G [q_1^T(x) \; q_2^T(x)] \begin{bmatrix} w_1(x) \\ w_2(x) \end{bmatrix} dx$$

(2.1)

We set

$$w(x) := \begin{bmatrix} w_1(x) \\ w_2(x) \end{bmatrix}, \quad K(x) := \begin{bmatrix} K_{11}(x) & K_{12}(x) \\ K_{21}(x) & K_{22}(x) \end{bmatrix}, \quad L(x,y) := \begin{bmatrix} L_{11}(x,y) & L_{12}(x,y) \\ L_{21}(x,y) & L_{22}(x,y) \end{bmatrix},$$

$$q(x) := \begin{bmatrix} q_1(x) \\ q_2(x) \end{bmatrix}$$

(2.2)

We postulate:

[A]. The matrix-valued functions (with real entries) $K(x)$, $L(x,y)$ are bounded measurable on $G$, $G \times G$, respectively.

[B]. The matrix-valued functions $K(x)$, $L(x,y)$ are symmetric, in the sense that
$K_{11}^T(x) = K_{11}(x)$, $K_{22}^T(x) = K_{22}(x)$, $K_{12}^T(x) = K_{21}(x)$,
$L_{11}^T(x,y) = L_{11}(y,x)$, $L_{22}^T(x,y) = L_{22}(y,x)$, $L_{12}^T(x,y) = L_{21}(y,x)$.

(This entails no loss of generality: arbitrary (i.e. not necessarily symmetric) matrices $K_{ij}(x)$, $L_{ij}(x,y)$ can be replaced by their symmetrizations
$\tilde{K}_{ij}(x) := \tfrac{1}{2}[K_{ij}(x) + K_{ji}^T(x)]$, $\tilde{L}_{ij}(x,y) := \tfrac{1}{2}[L_{ij}(x,y) + L_{ji}^T(y,x)]$ without affecting the values of the corresponding quadratic integral forms.)

[C]. Each of the matrices $K(x)$, $K_{11}(x)$, $K_{22}(x)$ is nonsingular, with an inverse that is also bounded measurable on $G$.

[D]. The matrix-valued functions $(K_{11}(x), L_{11}(x,y))$ are <u>jointly positive definite</u> in the sense that, for every non-zero $w_1 \in L^2(G, R^m)$, $\int_G w_1^T(x) K_{11}(x) w_1(x)\, dx + \int_G \int_G w_1^T(x) L_{11}(x,y) w_1(y)\, dx\, dy > 0$.

The matrix-valued functions $(K_{22}(x), L_{22}(x,y))$ are <u>jointly negative definite</u> in the sense that, for every non-zero $w_2 \in L^2(G, R^n)$, $\int_G w_2^T(x) K_{22}(x) w_2(x)\, dx + \int_G \int_G w_2^T(x) L_{22}(x,y) w_2(y)\, dx\, dy < 0$.

We shall prove:

<u>Proposition 2.1.</u> For every $w_2 \in L^2(G, R^n)$, there exists a unique minimizer of $E$ with respect to $w_1$. (We shall denote that minimizer by $(w_{1*} | w_2)$.)

<u>Proof:</u> For arbitrary but, in the context of the present proof, fixed $w_2 \in L^2(G, R^n)$, we set

$$E_1(w_1) := \tfrac{1}{2} \int_G w_1^T(x) K_{11}(x) w_1(x)\, dx + \tfrac{1}{2} \int_G \int_G w_1^T(x) L_{11}(x,y) w_1(y)\, dx\, dy + \\ + \int_G w_1^T(x) K_{12}(x) w_2(x)\, dx + \int_G \int_G w_1^T(x) L_{12}(x,y) w_2(y)\, dx\, dy + \int_G q_1^T(x) w_1(x)\, dx \quad (2.3)$$

The first variation of $E_1$ is

$$\delta E_1 = \int_G w_1^T(x) K_{11}(x) \delta w_1(x)\, dx + \int_G \int_G w_1^T(y) L_{11}(y,x) \delta w_1(x)\, dy\, dx + \\ + \int_G w_2^T(x) K_{21}(x) \delta w_1(x)\, dx + \int_G \int_G w_2^T(y) L_{21}(y,x) \delta w_1(x)\, dy\, dx + \int_G q_1^T(x) \delta w_1(x)\, dx \quad (2.4)$$

For a solution $\tilde{w}_1$ that annihilates $\delta E_1$, the second variation of $E_1$ is

$$\delta^2 E_1 = \int_G \delta w_1^T(x) K_{11}(x) \delta w_1(x)\, dx + \int_G \int_G \delta w_1^T(y) L_{11}(y,x) \delta w_1(x)\, dy\, dx \quad (2.5)$$

thus

$$\delta^2 E_1 > 0 \text{ if } \delta w_1 \text{ is not a.e. 0 in } G \quad (2.6)$$

by the condition of joint positive definiteness of the pair of kernels $(K_{11}, L_{11})$.

It remains to show the existence and uniqueness of a $\tilde{w}_1$ that annihilates $\delta E_1$, i.e. solves

$$K_{11}(x) w_1(x) + \int_G L_{11}(x,y) w_1(y)\, dy + q_1(x) + K_{12}(x) w_2(x) + \int_G L_{12}(x,y) w_2(y)\, dy = 0 \quad (2.7)$$

Since $K_{11}(x)$ is nonsingular, eq. (2.7) is a linear Fredholm integral equation of the second kind in the unknown function $w_1$ and the Fredholm alternative applies: either () has a solution, or the corresponding homogeneous equation has a nontrivial solution, i.e. for some nonzero $w_1$ we have

$$K_{11}(x)w_1(x) + \int_G L_{11}(x,y)w_1(y)\,dy = 0 \tag{2.8}$$

Pre-multiplication of (2.8) by $w_1^T(x)$ and integration (with $x$ as the variable of integration) of the resulting equation over $G$ yields a contradiction to the joint positive-definiteness of $(K_{11}, L_{11})$. Thus (2.7) is solvable. The same type of positive-definiteness argument shows that the solution is unique. ///

<u>Notation</u>: we shall denote the solution obtained in the above proposition by $(w_{1*} | w_2)$.

<u>Corollary 2.2.</u> By the same token, for every $w_1 \in L^2(G, R^M)$, there is a unique maximizer $(w_2^* | w_1)$ of $E$ with respect to $w_2$. ///

<u>Theorem 2.1.</u> The mapping from $L^2(G, R^m \times R^n)$ into itself, defined by

$$\begin{bmatrix} w_1 \\ w_2 \end{bmatrix} \mapsto \begin{bmatrix} (w_{1*} | w_2) \\ (w_2^* | w_1) \end{bmatrix} \tag{2.9}$$

has a unique fixed point.

<u>Notation</u>: We shall denote that unique fixed point by $\begin{bmatrix} w_{1**} \\ w_2^{**} \end{bmatrix}$.

<u>Proof</u>: We consider the system

$$\Theta_1 E = 0 \,;\quad \Theta_2 E = 0 \tag{2.10}$$

where

$$\delta E = \int_G [\Theta_1 E(w_1, w_2; x)\delta w_1(x) + \Theta_2 E(w_1, w_2; x)\delta w_2(x)]\,dx \tag{2.11}$$

The system (2.10) becomes

$$K(x)w(x) + \int_G L(x,y)w(y)\,dy + q(x) = 0 \tag{2.12}$$

and *in extenso*

$$\begin{aligned}
K_{11}(x)w_1(x) + K_{12}(x)w_2(x) + \int_G [L_{11}(x,y)w_1(y) + L_{12}(x,y)w_2(y)]\,dy + q_1(x) = 0\,; \\
K_{21}(x)w_1(x) + K_{22}(x)w_2(x) + \int_G [L_{21}(x,y)w_1(y) + L_{22}(x,y)w_2(y)]\,dy + q_2(x) = 0
\end{aligned} \tag{2.13}$$

By our assumption on non-singularity of $K(x)$, eq. (2.12) is a linear Fredholm integral equation of the second kind, and we establish existence of solution via the Fredholm alternative. It suffices to show that the associated homogeneous equation

$$K(x)w(x) + \int_G L(x,y)w(y)\,dy = 0 \tag{2.14}$$

cannot have a non-zero solution. If it does, we pre-multiply (2.14) by $\hat{w}^T(x)$ and integrate over $G$ with respect to $x$, where

$$\hat{w}(x) := \begin{bmatrix} w_1(x) \\ -w_2(x) \end{bmatrix} \tag{2.15}$$

thus obtaining

$$\int_G [w_1^T(x)K_{11}(x)w_1(x) - w_2^T(x)K_{22}(x)w_2(x)]\,dx +$$
$$+ \int_G \int_G [w_1^T(x)L_{11}(x,y)w_1(y) - w_2^T(x)L_{22}(x,y)w_2(y)]\,dy\,dx = 0 \tag{2.16}$$

which contradicts the conjunction of the joint positive-definiteness of $(K_{11}, L_{11})$ and the joint negative-definiteness of $(K_{22}, L_{22})$. Thus () has a solution which we denote by $\begin{bmatrix} w_{1**} \\ w_2^{**} \end{bmatrix}$. It is plain, from the extensive form () of (), that $w_{1**} = (w_{1*} \mid w_2^{**})$ and $w_2^{**} = (w_2^* \mid w_{1**})$, which proves the existence and the uniqueness of the wanted fixed point. ///

### 3. Jointly definite pairs of kernels.

We examine the quadratic form on $L^2(G, R^n)$ defined by a pair of kernels $K(x)$ and $L(x,y)$:

$$\Phi(w) := \int_G w^T(x)K(x)w(x)\,dx + \int_G \int_G w^T(x)L(x,y)w(y)\,dx\,dy \tag{3.1}$$

We assume the same conditions of measurability, boundedness, and symmetry, as in section 2.

The form $\Phi$ can be represented as follows:

$$\Phi(w) = \tfrac{1}{2}\int_G \int_G \tilde{w}^T(x,y) M(x,y)\tilde{w}(x,y)\,dx\,dy \tag{3.2}$$

where

$$\tilde{w}(x,y) := \begin{bmatrix} w(x) \\ w(y) \end{bmatrix}, \quad M(x,y) := \begin{bmatrix} |G|^{-1} K(x) & L(x,y) \\ L(x,y) & |G|^{-1} K(y) \end{bmatrix}; \quad |G| := \int_G dx \tag{3.3}$$

The <u>proof</u> of (3.2) is a direct calculation of substituting (3.3) into (3.2).

Consequently the positive-definiteness of $M(x,y)$, uniformly for $(x,y) \in G \times G$, is a sufficient condition for the joint positive-definiteness of $(K,L)$.

If we assume, in addition to all our other assumptions, that $K(x)$ is positive definite, uniformly for $x \in G$, then the nonnegative-definiteness of the form

$$\Psi(w) := \int_G \int_G w^T(x) L(x,y) w(y)\,dx\,dy \tag{3.4}$$

is a sufficient condition for the joint positive-definiteness of $(K,L)$.

When the kernel $L$ is continuous on $\overline{G \times G}$, by the same type of argument as that used in the well-known theory associated with the name of Mercer [M], the nonnegative-definiteness of $L$ is equivalent to the condition that, for all finite collections of points $\{x_i : 1 \leq i \leq N\}$ in $G \times G$, and $n$ – dimensional real vectors $\{a_i, 1 \leq i \leq N\}$, we should have

$$\sum_{i=1}^{N}\sum_{j=1}^{N} a_i^T L(x_i, x_j) a_j \geq 0 \tag{3.5}$$

When $L$ has the representation

$$L(x,y) = \sum_{(k,\ell)\in \boldsymbol{n}^2} \lambda_{k\ell}\omega_k(x)\omega_\ell^T(y) \tag{3.6}$$

where $\boldsymbol{n}$ is a countably infinite set of indices and $\{\omega_k : k \in \boldsymbol{n}\}$ is a complete orthonormal system in $L^2(G, R^n)$, then the nonnegative definiteness of $L$ is tantamount to the nonnegative-definiteness of the infinite-dimensional matrix $\Lambda := [\lambda_{k\ell}]_{(k,\ell)\in \boldsymbol{n}^2}$ in $\ell^2(\boldsymbol{n}, R)$. This is verified by expanding a generic element $w$ of $L^2(G, R^n)$ as

$$w(x) = \sum_{k\in \boldsymbol{n}} c_k \omega_k(x) \tag{3.7}$$

after which it follows that

$$\int_G \int_G w^T(x) L(x,y) w(y)\, dx\, dy = \sum_{(k,\ell)\in \boldsymbol{n}^2} \lambda_{k\ell} c_k c_\ell \tag{3.8}$$

## 4. Application to linear-quadratic Volterra integral games.

We consider a linear Volterra integral game with state equation

$$y(t) = y_0(t) + \int_{t_0}^{t} [A(t,s)y(s) + B(t,s)u(s) + C(t,s)v(s)]ds \qquad (4.1)$$

The state $y$ takes values in $p$ – dimensional Euclidean space, and the two controls, $u$ (minimizing) and $v$ (maximizing), in $m$ – and $n$ – dimensional Euclidean spaces, respectively.

The <u>time horizon</u> of our problem will be an interval $[t_0, t_1]$.

It is clear that our model covers also the case of linear differential games.

By using a <u>resolvent kernel</u> $S(t,s)$ corresponding to direct kernel $A(t,s)$, the solution of eq. (4.1) takes the form

$$y(t) = y_1(t) + \int_{t_0}^{t} [B_1(t,s)u(s) + C_1(t,s)v(s)]ds \qquad (4.2)$$

where

$$y_1(t) = y_0(t) + \int_{t_0}^{t} S(t,\sigma) y_0(\sigma) d\sigma ;$$

$$B_1(t,s) = B(t,s) + \int_{s}^{t} S(t,\sigma) B(\sigma,s) d\sigma ; \qquad (4.3)$$

$$C_1(t,s) = C(t,s) + \int_{s}^{t} S(t,\sigma) C(\sigma,s) d\sigma$$

It is convenient to take (4.2) as the equation of state dynamics and as the starting point of all further calculations.

By using the representation (4.2), it is clear that a quadratic functional in $y$, $u$, and $v$ can be transformed to a quadratic functional in $u$ and $v$ of the type discussed in section 2 of the present paper. Then all the mathematical questions boil down to finding conditions that guarantee the appropriate positive – and negative – definiteness conditions that imply existence and uniqueness of a saddle-point control.

For simplicity we will consider only a particular case of quadratic performance functional, that still provides enough structure to create a meaningful game-theoretic problem. We take a performance functional

$$J := \tfrac{1}{2} y^T(t_1) P_0 y(t_1) + \tfrac{1}{2} \int_{t_0}^{t_1} [y^T(t)P_1(t)y(t) + u^T(t)Q_1(t)u(t) + v^T(t)R_1(t)v(t)]dt +$$
$$+ \tfrac{1}{2} \int_{t_0}^{t_1} \int_{t_0}^{t_1} [y^T(t)P_2(t,s)y(s) + u^T(t)Q_2(t,s)u(s) + v^T(t)R_2(t,s)v(s)]ds\,dt \qquad (4.4)$$

We will assume, without expending or repeating the definitions, that all kernels that appear in eq. (4.4) have the appropriate symmetry properties as detailed in section 2.

By using the representation (4.2), we can represent $J$ as

$$J = \tfrac{1}{2}\int_{t_0}^{t_1} [u^T(t)\ v^T(t)] \begin{bmatrix} Q_1(t) & 0 \\ 0 & R_1(t) \end{bmatrix} \begin{bmatrix} u(t) \\ v(t) \end{bmatrix} dt +$$

$$+ \tfrac{1}{2}\int_{t_0}^{t_1}\int_{t_0}^{t_1} [u^T(t)\ v^T(t)] \begin{bmatrix} L_{11}(t,s) & L_{12}(t,s) \\ L_{21}(t,s) & L_{22}(t,s) \end{bmatrix} \begin{bmatrix} u(s) \\ v(s) \end{bmatrix} ds\, dt + \qquad (4.5)$$

$$+ \int_{t_0}^{t_1} [q_1^T(t)u(t) + q_2^T(t)v(t)]\, dt + \Omega$$

The functional $\Omega$ involves terms independent of the controls $u$ and $v$, and therefore does not affect the solution of the game-theoretic problem; for this reason, we shall not give the explicit calculation of $\Omega$.

The kernels in (4.5) are expressed in terms of the kernels in (4.4) and (4.2) as follows:

$$L_{11}(t,s) = B_1^T(t_1,t) P_0 B_1(t_1,s) + \int_{\max(t,s)}^{t_1} B_1^T(\sigma,t) P_1(\sigma) B_1(\sigma,s)\, d\sigma +$$

$$+ \int_t^{t_1}\int_s^{s_1} B_1^T(\tau,t) P_2(\tau,\sigma) B_1(\sigma,s)\, d\sigma\, d\tau + Q_2(t,s);$$

$$L_{12}(t,s) = B_1^T(t_1,t) P_0 C_1(t_1,s) + \int_{\max(t,s)}^{t_1} B_1^T(\sigma,t) P_1(\sigma) C_1(\sigma,s)\, d\sigma +$$

$$+ \int_t^{t_1}\int_s^{t_1} B_1^T(\tau,t) P_2(\tau,\sigma) C_1(\sigma,s)\, d\sigma\, d\tau;$$

$$L_{22}(t,s) = C_1^T(t_1,t) P_0 C_1(t_1,s) + \int_{\max(t,s)}^{t_1} C_1^T(\sigma,t) P_1(\sigma) C_1(\sigma,s)\, d\sigma +$$

$$+ \int_t^{t_1}\int_s^{s_1} C_1^T(\tau,t) P_2(\tau,\sigma) C_1(\sigma,s)\, d\sigma\, d\tau + R_2(t,s);$$

$$q_1^T(t) = y_1^T(t_1) P_0 B_1(t_1,t) + \int_t^{t_1} y_1^T(\sigma) P_1(\sigma) B_1(\sigma,t)\, d\sigma + \int_{t_0}^{t_1}\int_t^{t_1} y_1^T(\tau) P_2(\tau,\sigma) B_1(\sigma,t)\, d\sigma\, d\tau;$$

$$q_2^T(t) = y_1^T(t_1) P_0 C_1(t_1,t) + \int_t^{t_1} y_1^T(\sigma) P_1(\sigma) C_1(\sigma,t)\, d\sigma + \int_{t_0}^{t_1}\int_t^{t_1} y_1^T(\tau) P_2(\tau,\sigma) C_1(\sigma,t)\, d\sigma\, d\tau$$

(4.6)

It can be readily verified that, if the kernels $P_1(\sigma), P_2(\sigma,\tau)$ are nonnegative-definite (in the sense detailed in section 2), then $L_{11}$ is nonnegative-definite; likewise, if $P_1(\sigma), P_2(\sigma,\tau)$ are nonpositive-definite, then $L_{22}$ is nonpositive-definite. The situation, however, is complicated by the fact that the same kernels $P_1(\sigma), P_2(\sigma,\tau)$ appear in the expressions for both, $L_{11}$ and $L_{22}$. Consequently, we do not rely on $P_1(\sigma), P_2(\sigma,\tau)$, but we examine the kernels $Q_1, R_1, Q_2, R_2$. For the kernels $Q_1, R_1$, the <u>coercivity</u> constants $\alpha(Q_1), \alpha(R_1), \alpha(Q_1, R_1)$ are nonnegative numbers defined by the requirements that, for every $w \in L^2(0,T; R^m)$, we have

$$\int_{t_0}^{t_1} w^T(t) Q_1(t) w(t) \, dt \geq \alpha(Q_1) \| w \|^2_{L^2(0,T; R^m)} \, ;$$

$$\int_{t_0}^{t_1} \int_{t_0}^{t_1} w^T(t) R_1(t,s) w(s) \, ds \, dt \geq \alpha(R_1) \| w \|^2_{L^2(0,T; R^m)} \, ;$$ (4.7)

$$\int_{t_0}^{t_1} w^T(t) Q_1(t) w(t) \, dt + \int_{t_0}^{t_1} \int_{t_0}^{t_1} w^T(t) R_1(t,s) w(s) \, ds \, dt \geq \alpha(Q_1, R_1) \| w \|^2_{L^2(0,T; R^m)}$$

Of course, these coercivity constants are meaningful under the appropriate conditions of nonnegative-definiteness or positive-definiteness of $Q_1, R_1, (Q_1, R_1)$, respectively. Likewise, under appropriate nonpositive-definiteness conditions or negative-definiteness conditions, the <u>accretivity</u> constants $\beta(Q_2), \beta(R_2), \beta(Q_2, R_2)$, defined by

$$\int_{t_0}^{t_1} w^T(t) Q_2(t) w(t) \, dt \leq -\beta(Q_2) \| w \|^2_{L^2(0,T; R^n)} \, ;$$

$$\int_{t_0}^{t_1} \int_{t_0}^{t_1} w^T(t) R_2(t,s) w(s) \, ds \, dt \leq -\beta(R_2) \| w \|^2_{L^2(0,T; R^n)} \, ;$$ (4.8)

$$\int_{t_0}^{t_1} w^T(t) Q_2(t) w(t) \, dt + \int_{t_0}^{t_1} \int_{t_0}^{t_1} w^T(t) R_2(t,s) w(s) \, ds \, dt \leq -\beta(Q_2, R_2) \| w \|^2_{L^2(0,T; R^n)}$$

Consequently, sufficient conditions for $(Q_1, R_1, L_{11}. L_{22})$ to fulfill the conditions of section 2 (with $Q_1, R_1$ in the role of $K_{11}, K_{22}$) are:

(I). $\alpha(Q_1) > 0$, $\beta(R_1) > 0$.
(II). $\alpha(Q_1, Q_2)$, $\beta(R_1, R_2)$ are sufficiently large.

Of course, "sufficiently large" means large enough to make $L_{11}$ nonnegative-definite and $L_{22}$ nonpositive-definite. The detailed calculations, of how large these constants should be, are straightforward, and we do not include them here.

5. Volterra integral games that are linear-quadratic in the controls (LQC).

We consider a game with state equation that is linear in the controls, with performance functional that is quadratic in the controls; the state dynamics and the performance functional are generally nonlinear in the state.
The controlled integral equation is

$$y(t) = y_0(t) + \int_{t_0}^{t} [f_0(t,s,y(s)) + F_1(t,s,y(s))u(s) + F_2(t,s,y(s))v(s)]\,ds \tag{5.1}$$

and the performance functional is

$$J := \int_{t_0}^{t_1} \{g_0(t,y(t)) + g_1(t,y(t))u(t) + g_2(t,y(t))v(t) + \\ + \tfrac{1}{2}u^T(t)G_{11}(t)u(t) + u^T(t)G_{12}(t)v(t) + \tfrac{1}{2}v^T(t)G_{22}(t)v(t)\}\,dt \tag{5.2}$$

We set

$$G_{21}(t,y(t)) := G_{12}^T(t,y(t)) \tag{5.3}$$

As in the previous sections, $u$ is the minimizing control, and $v$ the maximizing control. $\psi$ will stand for the co-state; its values are co-vectors (row vectors) of the same dimension as the state $y$. The Hamiltonian of the problem is

$$H(t,y,u,v,\psi(\cdot)) := g_0(t,y) + g_1(t,y)u + g_2(t,y)v + \\ + \tfrac{1}{2}u^T G_{11}(t,y)u + u^T G_{12}(t,y)v + \tfrac{1}{2}v^T G_{22}(t,y)v + \\ + \int_{t}^{t_1} \psi(s)[f_0(s,t,y) + F_1(s,t,y)u + F_2(s,t,y)v]\,ds \tag{5.4}$$

For each control function $v$, we denote by $(u_* | v)$ a minimizer of $J$ over u, and by $v^{**}$ a maximizer of $J((u_* | v), v)$ over $v$. We set $u_{**} := (u_* | v^{**})$. Similarly $(v^* | u)$ s a maximizer of $J(u,v)$ over $v$, and $\bar{u}_{**}$ is a minimizer of $J(u,(v^* | u))$ over $u$, and we set $\bar{v}^{**} := (v^* | \bar{u}_{**})$. A pair of controls $(u_{**}, v^{**})$ gives the <u>lower value</u> of the game, and a pair $(\bar{u}_{**}, \bar{v}^{**})$ gives the upper value of the game.

We will carry out the construction of $(u_{**}, v^{**})$. We start with constructing $(u_* | v)$.

The first construction will rely on the Hamiltonian equation for the costate and an extremum principle for the control $u$. For each control function $v$, this is an optimal control problem for a Volterra integral system, for which necessary conditions for optimality are well-known. The two equations we need are

$$\psi(t) = \nabla_y H(t, y(t), u(t), v(t), \psi(\cdot)) ;$$
$$\nabla_u H(t, y(t), u(t), v(t), \psi(\cdot)) = 0 \tag{5.5}$$

From the extremum principle for $u$ (i.e. the second equation of (5.5)), we find

$$(u_* \mid v)(t) = -G_{11}^{-1}(t, y(t)) \left[ G_{12}(t, y(t)) v(t) + g_1^T(t, y(t)) + \int_t^{t_1} F_1^T(s, t, y(t)) \psi^T(s) \, ds \right] \tag{5.6}$$

Substitution of (5.6) into the Hamiltonian equation (the first of (5.5)) yields the <u>restricted Hamiltonian</u>

$$H_*(t, y, v, \psi(\cdot)) \equiv H(t, y, (u_* \mid v), v, \psi(\cdot)) =$$

$$= g_0(t, y) - g_1(t, y) G_{11}^{-1}(t, y) \left[ G_{12}(t, y) v + \int_t^{t_1} F_1^T(s, t, y) \psi^T(s) \, ds \right] +$$

$$+ (-1) \left[ g_1(t, y) + v^T G_{21}(t, y) + \int_t^{t_1} \psi(s) F_1(s, t, y) \, ds \right] G_{11}^{-1}(t, y) G_{12}(t, y) v +$$

$$+ \tfrac{1}{2} v^T G_{22}(t, y) v + \tfrac{1}{2} \left[ g_1(t, y) + v^T G_{21}(t, y) + \int_t^{t_1} \psi(s_1) F_1(s_1, t, y) \, ds_1 \right] \cdot$$

$$\cdot G_{11}^{-1}(t, y) \left[ g_1^T(t, y) + G_{12}(t, y) v + \int_t^{t_1} f_1^T(s_2, t, y) \psi^T(s_2) \, ds_2 \right] +$$

$$+ \int_t^{t_1} \psi(s) [f_0(s, t, y) + F_2(s, t, y) v] \, ds +$$

$$+ (-1) \int_t^{t_1} \psi(s) \left\{ F_1(s, t, y) G_{11}^{-1}(t, y) \left[ G_{12}(t, y) v + g_1^T(t, y) + \int_t^{t_1} F_1^T(\sigma, t, y) \psi^T(\sigma) \, d\sigma \right] \right\} ds \tag{5.7}$$

This restricted Hamiltonian is a quadratic functional in $v$ and consequently a critical point is

$$\underline{v}^{**}(t) = [G_{22}(t, y(t)) - G_{21}(t, y(t)) G_{11}^{-1}(t, y(t)) G_{12}(t, y(t))]^{-1} \cdot$$

$$\cdot \left\{ G_{21}(t, y(t)) g_1^T(t, y(t)) - g_2^T(t, y(t)) + \right. \tag{5.8}$$

$$\left. + \int_t^{t_1} [G_{21}(t, y(t)) G_{11}^{-1}(t, y(t)) F_1^T(s, t, y(t)) - F_2^T(s, t, y(t))] \psi^T(s) \, ds \right\}$$

This critical point will be unique and a maximizer of the restricted Hamiltonian $H_*$ under certain conditions to be specified later.

For the purpose of simplifying the notation, we shall make the convention that, unless otherwise explicitly specified, the variables inside each matrix-valued function will be $t$ and $y(t)$.

We set

$$G_3 := G_{22} - G_{21}G_{11}^{-1}G_{12} ;$$
$$F_3^T(s,t,y) := G_{21}G_{11}^{-1}F_1^T(s,t,y) - F_2^T(s,t,y) ; \qquad (5.9)$$
$$g_3^T := G_{21}g_1^T - g_2^T ;$$
$$g_4^T := G_3^{-1}g_3^T , \quad F_4^T(s,t,y) := G_3^{-1}F_3^T(s,t,y)$$

Then

$$\underline{v}^{**}(t) = g_4^T + \int_t^{t_1} F_4^T(s,t,y(t))\psi^T(s)\,ds \qquad (5.10)$$

Then we have

$$\underline{u}_{**}(t) \equiv (u_* \mid \underline{v}^{**})(t) = g_5^T + \int_t^{t_1} F_5^T(s,t,y(t))\psi^T(s)\,ds ;$$
$$g_5^T := -G_{11}^{-1}(G_{12}g_4^T + g_1^T) ; \qquad (5.11)$$
$$F_5^T(s,t,y) := -G_{11}^{-1}(G_{12}F_3^T(s,t,y) + F_1^T(s,t,y))$$

Sufficient conditions for the existence, uniqueness, and maximin property (relative to the Hamiltonian) of the solution $(\underline{u}_{**}, \underline{v}^{**})$ are that the kernels $G_{11}$, $G_{22}$ are, respectively, uniformly positive-definite and uniformly negative-definite; in that case, the kernel $G_3 = G_{22} - G_{12}^T G_{11}^{-1} G_{12}$ is also uniformly negative-definite.

Now we will substitute the expressions (5.10) and (5.11) for in the state dynamics and in the Hamiltonian equation for the co-state to obtain a system of integral equations for $(\underline{y}, \underline{\psi})$, the state and the co-state for the lower game. In the few equations below, we shall use simply $(y, \psi)$ for the state and the co-state of the lower game, instead of $(\underline{y}, \underline{\psi})$.

The equations for state and co-state of the lower game are

$$y(t) = y_0(t) + \int_{t_0}^t c_0(t,s,y(s))\,ds + \int_{t_0}^t \int_s^{t_1} C_1(t,\sigma,s,y(s))\psi^T(\sigma)\,d\sigma\,ds ;$$
$$\psi(t) = \nabla_{y(t)} \underline{H}(t,y(t),y(\cdot),\psi(\cdot)) ;$$
$$\underline{H}(t,y(t),y(\cdot),\psi(\cdot)) := h_0(t,y(t)) + \int_t^{t_1} h_1(s,t,y(t))\psi^T(s)\,ds + \qquad (5.12)$$
$$+ \tfrac{1}{2}\int_t^{t_1}\int_t^{t_1} \psi(s) H_2(s,\sigma,t,y(t),y(\sigma),y(s))\psi^T(\sigma)\,d\sigma\,ds$$

where

$$c_0(t,s,y(s)) = f_0(t,s,y(s)) + F_1(t,s,y(s))g_5^T(s,y(s)) + F_2(t,s,y(s))g_4^T(s,y(s));$$

$$C_1(t,\sigma,s,y(s)) = F_1(t,s,y(s))F_5^T(\sigma,s,y(s)) + F_2(t,s,y(s))F_4^T(\sigma,s,y(s));$$

$$h_0(t,y(t)) = g_1 g_5^T + g_2 g_4^T + g_5 G_{11} F_5^T + g_2 G_{11} F_4^T + g_1 F_5^T + g_2 F_4^T;$$

$$h_1(s,t,y(t)) = f_0^T + g_1 F_5^T + g_2 F_4^T + g_5 F_1^T + g_4 F_2^T + \tfrac{1}{2}[g_5 G_{11} g_5^T + g_5 G_{12} g_4^T + g_4 G_{21} g_5^T + g_4 G_{22} g_4^T];$$

$$H_2(\sigma,s,t,y(t),y(s),y(\sigma)) = \tfrac{1}{2}[F_5(s,t,\ldots)G_{11}(t,\ldots)F_5^T(\sigma,t,\ldots) + F_5 G_{12} F_4^T + F_4 G_{21} F_5^T + F_4 G_{22} F_4^T +$$

$$+ F_1(s,t,\ldots)F_5^T(s,\sigma,\ldots) + F_5(s,\sigma,\ldots)F_1^T(\sigma,t,\ldots) + F_2(s,t,\ldots)F_4^T(s,\sigma,\ldots) + F_4(s,\sigma,\ldots)F_2^T(\sigma,t,\ldots)]$$

(5.14)

Remark 5.1. Non-zero-sum Volterra games (with different performance functional for each of the players), and games with several players (instead of just two), can be treated by similar methods as in sections 4 and 5 of the present paper.
The case of two players in a non-zero-sum game, with two performance functionals,

$$y(t) = y_0(t) + \int_{t_0}^{t} f(t,s,y(s),u(s),v(s))\,ds\,;$$

$$J_i := F_i(T,y(T)) + \int_{t_0}^{T} G_i(t,y(t),u(t),v(t))\,dt\,;\ i=1,2$$

(5.15)

where $J_1$ is to be minimized and $J_2$ is to be maximized, and if we are interested in the lower value of the game, we would have two co-states and two Hamiltonians,

$$H_i = (\nabla_Y F_i(T,Y))f(T,t,y,u,v) + G_i(t,y,u,v) + \int_t^T \psi_i(s) f(s,t,y,u,v)\,ds \quad (i=1,2) \qquad (5.16)$$

and Hamiltonian equations

$$\psi_i(t) = \nabla_y H_i(t,y(t),u(t),v(t),\psi_i(\cdot)) \quad (i=1,2) \qquad (5.17)$$

When we are interested in the lower value of the game, we minimize the Hamiltonian $H_1$ with respect to $u$ (which is possible for LQ or LQC games), thus obtaining $(u_* | v)$, then we form the Hamiltonian $H_{2*}(t,y,v,\psi_2(\cdot)) = H_2(t,y,(u_*|v),v,\psi_2(\cdot))$ which is then to be maximized with respect to $v$, thus yielding $v^{**}$, and then $u_{**} = (u_* | v^{**})$.
The same process can be used for 3 or more players which alternate between minimizers and maximizers, say with controls $u$ (minimizer), $v$ (maximizer) and $w$ (minimizer), where the sequence of control decisions is $u, v, w$. From the 3 Hamiltonians and 3 co-states, we determine first $(u_* | v, w)$, then $(v^{**} | w)$, then $w_{***}$, then $v^{***} = (v^{**} | w_{***})$ and $u_{***} = (u_* | v^{***}, w_{***})$.   ///

6. Pursuit-evasion Volterra games.

It will facilitate our presentation to start with a more general Volterra game, and, after we develop some general equations, to work in more specific directions. We consider the game with state equation

$$y(t) = y_0(t) + \int_{t_0}^{t} f(t,s,y(s),u(s),v(s))\,ds \tag{6.1}$$

over a time-horizon $[t_0, t_1]$, where $t_1$ is not fixed but rather the game continues as long as a non-capture criterion

$$\Phi(t, y(t), u(t), v(t)) \neq 0 \tag{6.2}$$

is satisfied, and terminates as soon as the capture criterion

$$\Phi(t_1, y(t_1), u(t_1), v(t_1)) = 0 \tag{6.3}$$

is met.

Eq. (6.1) will be postulated to hold for $t \in [t_0, t_1)$. At the final (but not *a priori* fixed) time $t_1$, it is convenient to set

$$U := u(t_1),\ V := v(t_1) \tag{6.4}$$

and treat $U$ and $V$ as additional variables.

The state $y$ can represent the difference between two states, say $x$ and $z$, corresponding to the positions of the pursuer and the evader. The "capture", according to our definition, takes place when (6.3) is satisfied, and it amounts to driving the state and the controls in a desired manifold; in case this manifold consists if the single point 0, for the desired value of the state, then it becomes capture in the ordinary sense, as used in the classical literature on differential games and positional games.

For notational efficiency, we shall denote $y(t_1)$ also by $Y$. We will have occasion to use also $\dot{y}(t_1)$, and we set $W := \dot{y}(t_1)$.

As it will become clear in the work below, this is, in essence, a matter of notational convenience: the values $u(t_1)$, $v(t_1)$ would have to be taken into account as separate decision variables, regardless of the notation we may choose.

The performance functional $J$ will be taken in the form

$$J := F(t_1, y(t_1), u(t_1), v(t_1)) + \int_{t_0}^{t_1} G(t, y(t), u(t), v(t))\, dt \tag{6.5}$$

The inclusion, in the performance functional, of the values of the controls at the final time $t_1$ does not create complications, since, even if the capture criterion $\Phi$ was independent of the final values of the controls, the remaining ingredients of the problem would cause the final values of the controls to appear in the Hamiltonian equations; on the contrary, omission of the final values of the controls from $J$ would create conceptual complications, since then we would have values of the controls that appear in other locations in a set of necessary conditions but affect neither the state dynamics nor the performance functional.

After the final values of the controls have explicitly included in $J$, it creates no additional conceptual complications if they are also explicitly included in the capture criterion $\Phi$.

Towards the discovery of the relevant first-order conditions, we will use a penalty term, which for convenience we split into 3 different parts:

$$P_0 := \Psi \left[ y_0(t_1) - y(t_1) + \int_{t_0}^{t_1} f(t_1, t, y(t), u(t))\, dt \right];$$

$$P_1 := \omega \Phi(t_1, y(t_1), u(t_1), v(t_1)); \tag{6.6}$$

$$P_2 := \int_{t_0}^{t_1} \psi(t) \left[ y_0(t) - y(t) + \int_{t_0}^{t} f(t, s, y(s), u(s), v(s))\, ds \right] dt$$

Now the decision variables are $u, v, U, V$, and the state variables are $y, Y$. The time of capture, $t_1$, may be considered as a variable affecting the state dynamics.
The Lagrangian of our problem is

$$L := P_0 + P_1 + P_2 + J \tag{6.7}$$

The symbol $\tilde{\delta}$ will denote variations with respect to $y, Y, t_1$. We have

$$\tilde{\delta}L = \psi(t_1)\left[ (y_0(t_1) - y(t_1)) + \int_{t_0}^{t_1} f(t_1, t, y(t), u(t), v(t))\, dt \right] \delta t_1 +$$

$$+ \Psi \left[ \frac{dy_0(t_1)}{dt_1} \delta t_1 - \delta Y + f(t_1, t_1, Y, U, V)\delta t_1 + \int_{t_0}^{t_1} \nabla_y f(t_1, t, y(t), u(t), v(t))\, \delta y(t)\, dt + \right.$$

$$\left. + \int_{t_0}^{t_1} f_{t_1}(t_1, t, y(t), u(t), v(t))\, \delta t_1\, dt \right] +$$

$$+ \int_{t_0}^{t_1} \psi(t)(-\delta y(t))\, dt + \int_{t_0}^{t_1}\int_{t}^{t_1} \psi(s) \nabla_y f(s, t, y(t), u(t), v(t))\, \delta y(t)\, ds\, dt +$$

$$+ \omega \left[ \frac{\partial \Phi(t_1, y(t_1),\ldots)}{\partial t_1} \delta t_1 + \left( \nabla_Y \Phi(t_1, y(t_1),\ldots) \right) \delta Y \right] +$$

$$+ \frac{\partial F(t_1, y(t_1),\ldots)}{\partial t_1} \delta t_1 + \left( \nabla_Y F(t_1, y(t_1),\ldots) \right) \delta Y + G(t_1, Y, U, V)\, \delta t_1 + \int_{t_0}^{t_1} \nabla_y G(t, y(t), u(t), y(t))\, \delta y(t)\, dt$$

$$\tag{6.8}$$

The Hamiltonian is, therefore,

$$H \equiv H(t, t_1, y, u, v, Y, U, V, \psi(\bullet), \Psi, \omega) :=$$
$$= F(t_1, Y, U, V) + G(t, y, u, v) + \omega \Phi(t_1, Y, U, V) + \Psi f(t_1, t, y, u, v) + \tag{6.9}$$
$$+ \int_t^{t_1} \psi(s) f(s, t, y, u, v)\, ds$$

At this stage, we are in a position to state two (out of the) necessary conditions for optimal equilibrium:

the Hamiltonian equations for the co-states

$$[\psi(t)\ \ \omega] = [\nabla_y H\ \ \nabla_Y H]$$

and the variational capture condition at termination of the game.

$$y_0(t_1) - y(t_1) + \int_{t_0}^{t_1} f(t_1, t, y(t), u(t), v(t)) \, dt = 0 \; ;$$

$$\frac{dy_0(t_1)}{dt_1} + f(t_1, t_1, Y, U, V) + \int_{t_0}^{t_1} f_{t_1}(t_1, t, y(t), u(t), v(t)) \, dt = \dot{y}(t_1) \equiv W \tag{6.10}$$

Then the variational capture condition (involving the terms, out of $\tilde{\delta}L$, that contain a factor of $\delta t_1$) becomes

$$\Psi W + \omega \frac{\partial \Phi(t_1, y(t_1), \ldots)}{\partial t_1} + \frac{\partial F(t_1, \ldots)}{\partial t_1} + G(t_1, \ldots) = 0 \tag{6.11}$$

The two Hamiltonian equations can be rewritten in extensive form as

$$\psi(t) = \nabla_y G(t, y(t), u(t), v(t)) + \Psi \nabla_y f(t_1, t, y(t), u(t), v(t)) +$$
$$+ \int_t^{t_1} \psi(s) \nabla_y f(s, t, \ldots) \, ds \; ; \tag{6.12}$$
$$\Psi = \omega \nabla_Y \Phi(t_1, Y, U, V) + \nabla_Y F(t_1, Y, U, V)$$

It is plain that $\Psi$ and $\omega$ can be eliminated (under a certain transversality condition) from the full set of Hamiltonian equations and variational capture condition. Substitution of the second of (6.12) into (6.11) gives

$$\omega[(\nabla_Y \Phi(t_1, \ldots))W + \Phi_{t_1}] + (\nabla_Y F(t_1, \ldots))W + G(t_1, \ldots) = 0 \, . \tag{6.13}$$

The expression in square brackets, in (6.13) above, is the Eulerian derivative $\dfrac{D\Phi(t_1, \ldots)}{Dt_1}$.

Assuming that this derivative does not vanish (and this is a type of <u>transversality condition</u> at the moment of capture), we have

$$\omega = -\left[\frac{D\Phi(t_1, \ldots)}{Dt_1}\right]^{-1} \left[(\nabla_Y F(t_1, \ldots))W + G(t_1, \ldots)\right] \tag{6.14}$$

then substitution of that expression for $\omega$ into the second of (6.12) gives the value of $\Psi$:

$$\Psi = -\left[\frac{D\Phi(t_1, \ldots)}{Dt_1}\right]^{-1} \left[(\nabla_Y F(t_1, \ldots))W + G(t_1, \ldots)\right]\left(\nabla_Y \Phi(t_1, Y, U, V)\right) + \nabla_Y F(t_1, Y, U, V) \tag{6.15}$$

Example 6.1. A particular case is a <u>linear Volterra game</u> with quadratic performance functional and quadratic capture criterion. Because of the particularities of such a game, we shall rely only partly on the above treatment of a general Volterra game with pursuit-evasion ingredients, and we shall develop certain aspects of the solution of such a game directly.

We take the state equation

$$y(t) = y_0(t) + \int_{t_0}^{t} [A(t,s)y(s) + B(t,s)u(s) + C(t,s)v(s)]ds \qquad (6.16)$$

The capture criterion is taken as

$$\Phi(y(t_1)) \equiv \tfrac{1}{2} y^T(t_1) M\, y(t_1) = 0 \qquad (6.17)$$

which is, of course, equivalent to

$$y(t_1) \in \ker(M^{\tfrac{1}{2}}) \qquad (6.18)$$

Here, $M$ is a symmetric nonnegative-definite matrix (of dimensions compatible with the dimension of $y$), $M^{\tfrac{1}{2}}$ is the symmetric nonnegative-definite square root of $M$, and $\ker(M^{\tfrac{1}{2}})$ is the nullspace of $M^{\tfrac{1}{2}}$.
The performance functional $J$ is

$$J := \tfrac{1}{2} y^T(t_1) M_0 y(t_1) + \tfrac{1}{2} \int_{t_0}^{t_1} [y^T(t) M_1 y(t) + u^T(t) Q u(t) + v^T(t) R v(t)]dt \qquad (6.19)$$

We have (as can be easily verified)

$$\Psi = -\left(Y^T MW\right)^{-1} \left[ Y^T M_0 W + \tfrac{1}{2}(Y^T M_1 Y + U^T Q U + V^T R V) \right] \left[ Y^T (M + M_0) \right] \qquad (6.20)$$

The Hamiltonian equation for the costate $\psi$ is

$$\psi(t) = y^T(t) M_1 + \Psi\, A(t_1, t) + \int_{t}^{t_1} \psi(s) A(s,t)\, ds \qquad (6.21)$$

If $\Sigma(s,t)$ is a (matricial) resolvent kernel for the backward-time Volterra integral equation (6.21) (with direct kernel $A(s,t)$), then the costate can be represented as

$$\psi(t) = y^T(t) M_1 + \Psi A(t_1, t) + \int_{t}^{t_1} \left[ y^T(s) M_1 + \Psi A(t_1, s) \right] \Sigma(s,t)\, ds \qquad (6.22)$$

Then, by the methods explained in previous sections, the saddle-point solution $(u_{**}, v^{**})$ is found as

$$u_{**}^T(t) = -\left[\int_t^{t_1} \psi(s)B(s,t)\,ds + \Psi B(t_1,t))\right]Q^{-1} ;$$
$$v^{**T}(t) = -\left[\int_t^{t_1} \psi(s)C(s,t)\,ds + \Psi C(t_1,t))\right]R^{-1} \tag{6.23}$$

The costate $\psi$, as expressed via (6.22), is still a functional of the state trajectory $y$. (This is in addition to the costate's $\psi$ dependence on $\Psi$.) To remove the dependence on y, we use the representation (proved in section 3)

$$y(t) = y_1(t) + \int_{t_0}^t [B_1(t,s)u_{**}(s) + C_1(t,s)v^{**}(s)]\,ds \tag{6.24}$$

Substitution of (6.23) into (6.24) results in expressing the state $y$ in terms of $\psi$ and $\Psi$. The resulting expression for $y$, when substituted into (6.21), gives an integral equation for $\psi$ which contains also $\Psi$. We state the result of these calculations:

$$\psi(t) = y_1^T(t)M_1 + \Psi\left[A(t_1,t) - \int_{t_0}^t [B(t_1,s)Q^{-1}B_1^T(t,s) + C_1(t_1,s)R^{-1}C_1^T(t,s)]M_1\,ds\right] +$$
$$+ \int_t^{t_1} \psi(s)A(s,t)\,ds - \int_{t_0}^{t_1}\int_{t_0}^{\min(s,t)} \psi(s)[B(s,\sigma)Q^{-1}B_1^T(t,\sigma) + C(s,\sigma)R^{-1}C_1^T(t,\sigma)]M_1\,d\sigma\,ds \tag{6.25}$$

Thus the problem has been reduced to determining the final time $t_1$ and finding $Y, U, V, W$, as solutions to a fixed-point problem, since the obtained values of $(u_{**}, v^{**})$, and the associated trajectory $y(t)$, are required to satisfy

$$u_{**}(t_1) = U,\ v^{**}(t_1) = V,\ y(t_1) = Y,\ \dot{y}(t_1) = W \tag{6.26}$$

This translates to the system

$$U = -\Psi B(t_1,t_1)Q^{-1} ;$$
$$V = -\Psi C(t_1,t_1)R^{-1} ;$$
$$Y = y_1(t_1) + \int_{t_0}^{t_1} [B_1(t_1,s)u_{**}(s) + C_1(t_1,s)v^{**}(s)]\,ds ; \tag{6.27}$$
$$W = \dot{y}_1(t_1) + B_1(t_1,t_1)U + C_1(t_1,t_1)V + \int_{t_0}^{t_1} [B_{1,t_1}(t_1,s)u_{**}(s) + C_{1,t_1}(t_1,s)v^{**}(s)]\,ds$$

Also, the capture condition (6.18) needs to be satisfied.

The system consisting of (6.18), (6.20), (6.21), (6.23) – (6.25), and (6.27) is the full set of first-order necessary conditions for the pursuit – evasion Volterra game, as formulated in the present example 6.1. It is a simple observation, but perhaps worth stating explicitly, that the determination of the capture time $t_1$ is the hardest part of this problem, and the only nonlinear ingredient; the rest of the necessary conditions are linear in the state, costate, and control variables. ///